%% file: delignemilnor.tex
\documentclass[a4paper,11pt,leqno]{smfart}

\input{stylemoi2}

\usepackage[T1]{fontenc}

\title{Conjecture de Bloch et nombres de Milnor}

\begin{document}

\maketitle

\section*
{Introduction}
\subsection{} \label{not1}
Soient $S=\SP(A)$ un trait hensélien à corps résiduel algébriquement clos, de point fermé (resp. générique) $s$ (resp. $\eta$), $X$ un 
$S$-schéma plat, séparé de type fini, purement de dimension relative $n\in \NN$, et lisse en dehors d'un unique point fermé $x$ de la fibre spéciale $X_s$.
On suppose de plus $X$ régulier.
Soit 
\begin{equation}
\mu(X/S,\, x)=\mathrm{long}_{\OO_{X,x}}\ \sous{\Ext}^1(\Omega^1_{X/S},\OO_X)_x,
\end{equation}
le nombre de Milnor de $X$ en $x$ (\sga{7}{xvi}{1.2}).

Soient $\etb$ un point géométrique localisé en $\eta$, et $\ell$ un nombre premier inversible dans $\OO_S$ ; le complexe des cycles évanescents sur $X_s$, noté
$\ce(\FF_{\ell})$, est concentré en $x$ et à cohomologie constructible. Pour tout $\FF_{\ell}$-espace vectoriel de dimension finie $M$, muni d'une action
continue de $\gp(\eta,\etb)$, on note $\mathrm{dim\ tot}\,M$ l'entier $\dim_{\FF_{\ell}}(M) + \swan(M)$.
Dans \sga{7}{xvi}{1.9}, P. Deligne fait la conjecture suivante :

\begin{cnj}[« Deligne-Milnor »]\label{conjDM}
Sous les hypothèses précédentes, on a l'égalité :

\begin{equation}\label{formDM}
\mu(X/S,\, x)=(-1)^n \mathrm{dim\ tot}\, \ce(\FF_{\ell})_x.
\end{equation}
\end{cnj}

Cette conjecture est démontrée dans {\sl loc. cit.} dans les trois cas suivants : 
\begin{itemize}
\item $n=0$,
\item $X/S$ présente une singularité quadratique ordinaire en $x$,
\item $S$ est d'égale caractéristique. \end{itemize}

\subsection{} 
Plus généralement, cette conjecture a un sens dès que $\kappa(s)$ est parfait. Cependant, on ignore comment définir un second membre sans cette
hypothèse.
Notons que, d'après \cite{liannul}, $\ce(\FF_{\ell})$ est concentré en degré $n$, de sorte que le second membre de \ref{formDM}
est $\dim_{\FF_{\ell}}(\ce^n(\FF_{\ell})_x) + \swan(\ce^n(\FF_{\ell})_x)$.

\subsection{}\label{not2} Soient $S$ comme précédemment et $X$ un $S$-schéma régulier, plat, séparé de type fini, purement de dimension 
relative $n$, à fibre générique lisse.
Soit 
$$
\mathrm{Art}(X/S)=\mathrm{dim\ tot}\, \RG(X_s, \ce(\FF_{\ell})),
$$
le conducteur d'Artin. Si $X/S$ est propre, le complexe des cycles proches calcule la cohomologie de la fibre générique géométrique, et l'on a 
$$
\mathrm{Art}(X/S)=\chi(X_{\etb})-\chi(X_{s})+\swan\, \RG(X_{\etb},\FF_{\ell}).
$$
Dans \cite{blochconj}, S. Bloch définit une classe de Chern localisée 
$$
{c_{n+1}}^{X}_{X_s}(\Omega^1_{X/S})\in \mathrm{CH}_0(X_s),
$$
et fait la conjecture suivante :

\begin{cnj}[Bloch]\label{conjB}
Supposons de plus $X/S$ propre, on a :
\begin{equation}\label{formB}
\mathrm{Art}(X/S)=(-1)^n \deg {c_{n+1}}^{X}_{X_s}(\Omega^1_{X/S}).
\end{equation}
\end{cnj}

Cette conjecture est démontrée par S.~Bloch dans {\sl loc. cit.} pour $n=1$, et par K.~Kato et T.~Saito
si l'on suppose que $(X_s)_{\red}$ est un diviseur à croisement normaux (\cite{katosaitobloch}).

\subsection{} Dans~\cite{licotang}, L. Illusie définit les dérivés du foncteur non additif $\Lambda^{n+1}$. Si $Z_{X/S}\subset X_s$ désigne le lieu fermé
de non lissité de $f:X\ra S$, le complexe $\mathsf{L}\Lambda^{n+1}\Omega^1_{X/S}$ appartient à $\mathsf{D}^b_{Z_{X/S}}(X)_{\mathrm{parf}}$. 
La structure de schéma (non nécessairement réduit) sur l'espace $Z_{X/S}$ est explicitée dans la section suivante.
D'après T. Saito (\cite{tsaitoselfint},~2.3 et~\cite{tsaitoparite}, corrections), on a l'égalité 
\begin{equation}\label{egsaito}
\deg {c_{n+1}}^{X}_{X_s}(\Omega^1_{X/S})=\chi(X,\mathsf{L}\Lambda^{n+1}\Omega^1_{X/S})
,
\end{equation}
où $\chi$ désigne le composé $\K_{X_s}(X)\iso \K(X_s)\sr{f_{s*}}{\ra} \K(s)=\ZZ$.
Le terme de droite de~\ref{egsaito} nous permet de définir le second terme de~\ref{formB} 
en supposant seulement que $Z_{X/S}$ est propre sur $s$ et par là même d'émettre la conjecture suivante :

\begin{cnj}\label{conjB2}
Soient $S$ et $X$ comme dans~\ref{not2}. Supposons le lieu $Z_{X/S}$ de non lissité de $X/S$ propre sur $s$. On a :
$$
\mathrm{Art}(X/S)=(-1)^n \chi(X,\mathsf{L}\Lambda^{n+1}\Omega^1_{X/S}).
$$
\end{cnj}

Nous verrons plus bas que c'est une généralisation commune de~\ref{conjDM} et~\ref{conjB}.
Le résultat principal de cette note est le théorème suivant :

\begin{thm}\label{theoreme}
La conjecture~\ref{conjDM} se déduit de la conjecture~\ref{conjB}.
\end{thm}

On en tire le

\begin{crl}
La formule de Deligne-Milnor est valable en dimension relative $1$.
\end{crl}

\subsection{}
Nous vérifierons dans la section suivante que la conjecture~\ref{conjDM} est équivalente à la conjecture~\ref{conjB2} dans le cas où $Z=\{x\}$.
En particulier,~\ref{conjDM} est équivalent à~\ref{conjB} si $X/S$ est propre et présente une unique singularité dans la fibre spéciale.
Il serait intéressant de généraliser l'énoncé~\ref{theoreme} en une démonstration de l'implication~\ref{conjB} $\Rightarrow$~\ref{conjB2}.

\subsection{Remerciements}
Ils vont en premier lieu à Luc Illusie qui a eu la gentillesse de me faire part de cette question et
de sa note~\cite{linoteconjbloch}, dans laquelle on trouve le lien entre puissances extérieures dérivées 
et nombre de Milnor. Ses commentaires sur les versions
précédentes de ce texte en ont grandement amélioré la lisibilité. 
Enfin, je remercie Michel Raynaud à qui je dois la démonstration du~\ref{Michel} et Ofer Gabber qui a eu l'amabilité 
de relire ce texte et d'y déceler quelques erreurs.

\section{Nombre de Milnor et classe de Bloch}
Les hypothèses sont celles du~\ref{not2}.

\subsection{Description locale du lieu singulier}\label{lieusing}
Localement sur $X$ pour la topologie de Zariski, il existe un $S$-schéma lisse $P$ de dimension relative $n+1$ et
une $S$-immersion régulière $i:X\hra P$ (cf. par exemple~\cite{katosaitobloch},~\S 1). Pour simplifier l'écriture, nous notons encore $X$ 
un tel ouvert.
La suite exacte en traits pleins
\begin{equation}\label{sefond}
0\dashrightarrow \mc{N}_{X/P}\sr{d}{\ra} i^*\Om^1_{P/S}\ra \Om^1_{X/S}\ra 0
\end{equation}
est aussi exacte à gauche. En effet, le faisceau $\mc{N}_{X/P}$ est localement libre (de rang $1$) par hypothèse et l'exactitude
à gauche est valable en restriction à la fibre générique $X_{\eta}$, supposée lisse sur $\eta$.
L'image $\mc{J}_X$ du morphisme $\mc{N}_{X/P}\otimes_{\OO_X} ({i^*\Om^1_{P/S}})^{\vee}\ra \OO_X$ définit un sous-schéma fermé $$e:Z\hra X.$$ 
C'est l'idéal Jacobien $\mc{J}^n_{X/S}$ (cf.~\sga{7}{vi}{\S 5}) abstraitement défini comme l'idéal
de Fitting $\mathrm{Fitt}_n(\Om^1_{X/S})$. En particulier, il est indépendant du choix de $P$, ce qui résulte 
aussi de~\ref{isomutile}.
De la suite exacte~\ref{sefond}, on déduit la suite exacte :
\begin{equation}
\mc{N}_{X/P}\otimes_{\OO_X} i^*\Om^n_{P/S} \ra i^*\Om^{n+1}_{P/S}\ra \Om^{n+1}_{X/S}\ra 0.
\end{equation}
Par tensorisation avec le faisceau $(i^*{\Om^{n+1}_{P/S}})^{\vee}$, localement libre de rang $1$, on en déduit une suite exacte :
\beqt \label{seutile}
\mc{N}_{X/P} \otimes_{\OO_X} (i^*\Om^1_{P/S})^{\vee} \sr{d^{\vee}}{\ra} \OO_X \ra 
\Om^{n+1}_{X/S}\otimes_{\OO_X} (i^*{\Om^{n+1}_{P/S}})^{\vee}\ra 0.
\eeqt
Ainsi, on a un isomorphisme 
\beqt \label{isomutile}
\OO_Z=\OO_X/\mc{J}_X \iso \Om^{n+1}_{X/S}\otimes_{\OO_X} (i^*{\Om^{n+1}_{P/S}})^{\vee},
\eeqt

\subsection{Expression locale de $T^1_{X/S}=\sous{\Ext}^1(\Om^1_{X/S},\OO_X)$}\label{t1}
Plaçons-nous dans un ouvert affine convenable de $X$, comme dans le paragraphe précédent.
La résolution localement libre de $\Om^1_{X/S}$ permet de calculer le faisceau $T^1_{X/S}$.
En appliquant le foncteur $\sous{\Hom}(-,\OO_X)$ à~\ref{sefond}, on trouve la suite exacte :
$$
({i^*\Om^1_{P/S}})^{\vee}\ra \mc{N}_{X/P}^{\vee}\ra \sous{\Ext}^1(\Om^1_{X/S},\OO_X)\ra 0.
$$

Tensorisant~\ref{seutile} avec $\mc{N}_{X/P}^{\vee}$, on obtient un isomorphisme
\begin{align}\label{t1local}
   T^1_{X/S}= &\  \Om^{n+1}_{X/S}\otimes_{\OO_{X}} (i^*{\Om^{n+1}_{P/S}})^{\vee}\otimes_{\OO_{X}} \mc{N}_{X/P}^{\vee} \\
   = &\  \OO_Z \otimes_{\OO_X} \mc{N}_{X/P}^{\vee}.
\end{align}

Ainsi, $T^1_{X/S}$ a pour support $Z$, donc est de longueur finie sur $\OO_{X,x}$ si $Z_{\red}=\{x\}$.
En particulier, si $P=\aff^{n+1}_S$, et $0 \in X=V(f)$ est une singularité isolée, on retrouve la définition usuelle du nombre de Milnor donnée dans
\sga{7}{xvi}{\S 1} :
$$\mu(f)=\mathrm{long}_A\,A[t_1,\dots,t_{n+1}]_{(t_1,\dots,t_{n+1})}/\big(f,\frac{\partial f}{\partial t_1},\dots, \frac{\partial f}{\partial t_{n+1}}\big).$$

\subsection{Complexes de Koszul et dérivés des puissances extérieures}

\subsubsection{Rappels et notations}

Soient $R$ un anneau local régulier, $r$ un entier, et $u=(u_1,\dots,u_r):R^{r}\ra R$ une application $R$-linéaire.
Notons $e_1,\dots,e_r$ la base canonique de $R^r$ et 
$$\kos(u) : [0\ra \Lambda^{r}R^r\ra\cdots \ra \Lambda^{k+1}R^r\ra \Lambda^k R^r \ra \cdots\ra R^r\sr{u}{\ra} R\ra 0],$$
le complexe de Koszul usuel, où $R$ est placé en degré $0$ et $\Lambda^{k+1}R^r\ra \Lambda^k R^r$ est 
donné par 
$$x=e_{i_1}\wedge\cdots\wedge e_{i_{k+1}}\mapsto x\llcorner u = \sum_{j=1}^{k+1}(-1)^{j-1} u(e_{i_j})e_{i_1}\wedge\cdots \wedge \widehat{e_{i_j}}  \wedge e_{i_{k+1}}.$$
On a $\HH^0\big( \kos(u)\big)=R/u(R^r)$, et le morphisme canonique $\kos(u)\ra R/u(R^r)$ est un isomorphisme si, et seulement si, la suite $u$ est 
régulière ({\it complètement sécante} dans la terminologie de Bourbaki~\cite{bbkalg10}), c'est-à-dire $\mathrm{long}_R \big(R/u(R^r)\big)<+\infty$ si l'on suppose
de plus $r=\dim\, R$.

Dualement, pour tout morphisme $v:R\ra R^r$, on a le complexe 
$$
\kosdual(v) : [0\ra R\sr{v}{\ra} R^r\ra \cdots \ra \Lambda^k R^r \sr{v\wedge}{\ra} \Lambda^{k+1}R^r\ra \cdots \ra \Lambda^{r}R^r \ra 0],
$$
où $R$ est à nouveau placé en degré $0$.

Rappelons enfin la dualité de Koszul (cf. par exemple~\cite{eisenbud}, 17.15) : $\kos(u)^{\vee}\iso \kosdual(u^{\vee})$.


\subsubsection{} 

Soient $R$ et $r$ comme précédemment, et $\mc{C}_v :[R\sr{v}{\ra} R^r]$ un objet de $\mathsf{D}_{\mathrm{coh}}^-(R)$ (la catégorie
des complexes bornés supérieurement de $R$-modules dont les groupes de cohomologie sont de type fini), où $R$ est placé en degré~$-1$.

\begin{lmm}
Avec les notations précédentes, on a un isomorphisme dans $\mathsf{D}_{\mathrm{coh}}^-(R)$ :
$$
\mathsf{L}\Lambda^r(\mc{C}_v)=\kosdual(v)[r].
$$
\end{lmm}

\begin{proof}
D'après l'isomorphisme de Quillen (cf.~\cite{licotang},~{\sc i}.4.3.2), on a $\mathsf{L}\Lambda^r(\mc{C}_v)=\mathsf{L}\Gamma^r(\mc{C}_v[-1])[r]$, où
$\Gamma$ désigne le foncteur non additif « algèbre à puissances divisées » (les tenseurs symétriques).
Il est démontré dans {\sl loc. cit.},~{\sc viii}.2.1.2.1 que si $\mc{L}=[R\sr{a}{\ra} R^{r}]$, ses composantes étant placées
en degré $0$ et $1$, on a $\mathsf{L}\Gamma^r(\mc{L})=\kosdual(a)$. Le lemme en découle.
\end{proof}

On trouvera dans {\sl loc. cit.} des résultats plus généraux : cas des complexes à composantes plates, etc.

En particulier, il résulte de la dualité de Koszul que le morphisme canonique $\mathsf{L}\Lambda^r \mc{C}_v\ra \Lambda^r \HH^0(\mc{C}_v)$ est un isomorphisme
si $v$ est une suite régulière.

\subsection{Globalisation}\label{globalisation}
\subsubsection{}Sous les hypothèses de~\ref{not2}, on a :
\begin{enumerate}
\item l'idéal Jacobien $\mc{J}_X$ est l'annulateur du $\OO_X$-module $\Om^{n+1}_{X/S}$. Par la suite, nous noterons $e:Z_{X/S}=V(\mathrm{Ann}\,\Om^{n+1}_{X/S})\hra X$.
\item \label{2} le morphisme canonique $L_{X/S}\ra \Om^1_{X/S}$, où $L_{X/S}$ est le complexe cotangent défini dans~\cite{licotang},
est un isomorphisme (dans la catégorie dérivée adéquate).
\item on a un isomorphisme canonique 
\beqt \label{t1global}
T^1_{X/S}={e_*(e^*\Om^{n+1}_{X/S}})^{\vee}\otimes_{\OO_X} \mathrm{d\acute{e}t}(\Om^1_{X/S})
\eeqt
\item le morphisme canonique 
\beqt \mathsf{L}\Lambda^{n+1}\Omega^1_{X/S} \ra \Om^{n+1}_{X/S}
\eeqt 
est un isomorphisme si $Z_{X/S}$ est de dimension $0$.

\end{enumerate}
Le premier énoncé résulte de~\ref{isomutile}.
Le second est bien connu (cf.~\cite{katosaitobloch}, \S 1.5) et  justifie la définition que nous avons prise du faisceau $T^1_{X/S}$ dans 
le paragraphe précédent. L'isomorphisme~\ref{t1global} est une globalisation de~\ref{t1local}, que nous laissons au lecteur. 
(Voir~\cite{detdiv} pour la définition du déterminant d'un complexe parfait.) 
Rappelons cependant que localement, avec les notations de~\ref{sefond}, 
on a $\mathrm{d\acute{e}t}(\Om^1_{X/S})=i^*\Om^{n+1}_{P/S}\otimes {\mc{N}_{X/P}}^{\vee}$.
Le dernier point résulte, par localisation, des calculs locaux précédents : $\mathsf{L}\Lambda^{n+1}\Omega^1_{X/S}$ est acyclique hors du degré $0$.

\subsubsection{} Pour mémoire, signalons le résultat suivant.
Soit $\mc{L}$ le faisceau $\sous{\HH}^{-1} Le^*\Om^1_{X/S}$ considéré dans~\cite{katosaitobloch}. Sous les hypothèses de~\ref{lieusing},
le complexe $Le^*\Om^1_{X/S}$ est (localement) isomorphe au complexe $[e^*\mc{N}_{X/P}\sr{0}{\ra}(ie)^*\Om^1_{P/S}]$, si bien que 
$\mc{L}$ est (localement) isomorphe à $e^*\mc{N}_{X/P}$. 
Le faisceau inversible $\mc{L}^{\vee}$, localement isomorphe à $e^*\mc{N}_{X/P}^{\vee}$, est globalement isomorphe à $e^*T^1_{X/S}$.

\section{Compactification}

Le résultat principal est le suivant :

\begin{prp}\label{deformproj}
Sous les hypothèses de~\ref{conjDM}, et si l'on suppose de plus $S$ complet, 
il existe un $S$-schéma projectif et plat $Y$, purement de dimension relative $n$, lisse en dehors d'un unique
point fermé $y$ de la fibre spéciale $Y_s$, tel que les hensélisés stricts $X_{(x)}$ et $Y_{(y)}$ soient 
$S$-isomorphes.
\end{prp}

La démonstration fait l'objet des paragraphes~\ref{not3} à~\ref{prpfin}. Dans le cas de la dimension relative~$1$, un autre argument, dû à M. Raynaud,
est donné en~\ref{Michel}.

\subsection{}\label{not3} Comme $X$ est régulier et que toute $S$-immersion dans un $S$-schéma lisse est régulière, il existe un entier $r$ tel que
$(X,x)$ soit Zariski-localement isomorphe à $(V(\mathbf{f}),0)$, où $\mathbf{f}=(f_1,\dots,f_r)$ est une suite régulière de $A[t_1,\dots,t_{n+r}]=A[\mathbf{t}]$
et $\{0\}$ désigne l'origine de $\aff^{n+r}_s$. On suppose désormais $S$ complet, $X=V(\mathbf{f})$ et $x=0$.
Notons $\MM$ (resp. $\hat{\MM}$) l'idéal maximal en l'origine de $A[\mathbf{t}]$ (resp. du complété $A\[\mathbf{t}\]$)
et $\hat{\MM}_X$ (resp. $\MM_X$) l'image de $\hat{\MM}$ dans $R\fed \widehat{\OO_{X,x}}$ (resp. l'idéal maximal de $\OO_{X,x}$).
Enfin, supposons $\mu=\mu(X/S,\,x)>0$.

\begin{lmm}\label{deform} Sous les hypothèses de~\ref{not3}, il existe un entier $\lambda_{X,x}$ tel que 
pour toute suite $\mathbf{g}=(g_1,\dots,g_r)$ de $A[t_1,\dots,t_{n+r}]$, 
satisfaisant les $r$ relations de congruences $g_i-f_i\in \MM^{\lambda_{X,x}}$, les deux schémas strictement 
locaux $V(\mathbf{g})_{(0)}$ et $X_{(x)}$ soient $S$-isomorphes.
\end{lmm}

La démonstration se coupe en deux : une partie formelle (\ref{suffjet}), et une de descente
aux hensélisés.

\begin{lmm}[Suffisance des jets] \label{suffjet}
Pour tout $r$-uplet d'éléments $\mathbf{g}\in A\[\mathbf{t}\]^{r}$, 
satisfaisant $\mathbf{f}-\mathbf{g}\in (\hat{\MM}^{3\mu})^r$,
il existe $\mathbf{x}=(x_1,\dots,x_{n+r})$ dans $A\[\mathbf{t}\]^{n+r}$, tel que $\mathbf{x}\equiv \mathbf{t}\mod \hat{\MM}^2$ et 
$\mathbf{g}(\mathbf{x})=0$ dans $A\[\mathbf{t}\]/(\mathbf{f})$. 
En d'autres termes, il existe un $A$-isomorphisme {\it tangent à l'identité} : 
$A\[\mathbf{t}\]/(\mathbf{g})\iso A\[\mathbf{t}\]/(\mathbf{f})$,
défini par $t_i\mapsto x_i$.

\end{lmm}

Notons $\mathbf{f}^{\prime}$ l'application linéaire $A[\mathbf{t}]^{n+r} \ra A[\mathbf{t}]^r$, définie par les dérivées partielles 
$\frac{\partial f_i}{\partial t_j}$.
Par hypothèse, le nombre de Milnor 
$$\mu=\mathrm{long}_{R}\,R^r/\mathrm{Im}(\mathbf{f}'_R)\geq 1,$$
est fini (cf.~\ref{t1}), où $\mathbf{f}'_R$ désigne $\mathbf{f}'\otimes_{A[\mathbf{t}]} R$. (De même, nous noterons $\mathbf{f}_R$ l'image
de $\mathbf{f}$ dans $R^r$, etc.)
Ainsi, $({\hat{\MM}_X}^{\mu} R)^r\subset \mathbf{f}'_R R^{n+r}$.
On en déduit, pour tout $c\in \NN$, l'inclusion de sous-$R$-modules de $R^r$ :
$$({\hat{\MM}_X}^{\mu+c}R)^{r}\subset \mathbf{f}'_R \big(({\hat{\MM}_X}^{c}R)^{n+r}\big)\ (\star_c).$$
Remarquons que si $\mathbf{g}\in A\[\mathbf{t}\]^r$ satisfait les congruences $\mathbf{f}_R-\mathbf{g}_R\in ({\hat{\MM}_X}^{\mu+2})^{r}$, 
l'inclusion $(\star_0)$ est encore valable avec $\mathbf{g}_R$ à la place de $\mathbf{f}_R$.
Considérons $\mathbf{g}\in A\[\mathbf{t}\]^r$ tel que $\mathbf{g}-\mathbf{f} \in (\hat{\MM}^{3\mu})^{r}$ 
comme dans l'énoncé et tâchons de vérifier les conclusions de~\ref{suffjet}.
Soit $\varepsilon\in R^{n+r}$ ; la formule de Taylor pour $\mathbf{g}$ s'écrit :
$$
\mathbf{g}_R(\mathbf{t}_R+\varepsilon)=\mathbf{g}_R(\mathbf{t}_R)+\mathbf{g}'_R\cdot \varepsilon + \big( \text{termes quadratiques en\ } \varepsilon\big).
$$
En particulier, d'après $(\star_{2\mu})$, on peut trouver $\varepsilon_{[0]}\in ({\hat{\MM}_X}^{2\mu}R)^{n+r}$ tel que
$\mathbf{g}'_R\cdot \varepsilon_{[0]}=\mathbf{f}_R-\mathbf{g}_R$. La formule précédente montre qu'on a alors 
$\mathbf{g}_R(\mathbf{t}_R+\varepsilon_{[0]})=\mathbf{g}_R+\alpha_{[1]}$, où $\alpha_{[1]}\in ({\hat{\MM}_X}^{4\mu}R)^{n+r}$.
La fonction $\mathbf{g}_{[0]}\fed \mathbf{g}_R(\mathbf{t}_R+\varepsilon_{[0]})$ satisfait les inclusions $(\star_c)$, pour le même $\mu$,
car on a $\mathbf{g}_R-\mathbf{g}_{[0]}\in  ({\hat{\MM}_X}^{3\mu})^{r} \subset ({\hat{\MM}_X}^{\mu+2})^{r}$.
Par récurrence, on construit de proche en proche, une suite éléments $\varepsilon_{[i]}\in ({\hat{\MM}_X}^{(2^i+1)\mu}R)^{n+r}$, $i\geq 0$, telle que
$$\mathbf{g}_{[i]}\fed \mathbf{g}_R(\mathbf{t}_R+\varepsilon_{[0]}+\cdots+\varepsilon_{[i]})
\big(=\mathbf{g}_{[i-1]}(\mathbf{t}_R+\varepsilon_{[i]})\big)=\mathbf{f}_R+\alpha_{[i+1]},
$$
où $\alpha_{[i+1]}\in ({\hat{\MM}_X}^{(2^{i+1}+2)\mu}R)^{r}$.
L'anneau $R$ étant complet, on peut considérer
$$\varepsilon=\sum_{i=0}^{\infty} \varepsilon_{[i]}\in \hat{\MM}_X^2.$$ 
On a $\mathbf{g}_R(\mathbf{t}_R+\varepsilon)=\mathbf{f}_R(\mathbf{t}_R)$.
Donc, si $\tilde{\varepsilon}\in (\hat{\MM}^2)^{n+r}$ relève $\varepsilon$, alors $\mathbf{x}=\mathbf{t}+\tilde{\varepsilon}$ vérifie les 
conditions de~\ref{suffjet}.

\subsubsection*{Algébrisation}\label{algebris}
Montrons que l'entier $\lambda_{X,x}=3\mu$ de~\ref{suffjet} convient pour~\ref{deform}.
Soit $$B={A[t_1,\dots,t_{n+r}]_{\MM}}^{hs}/(f_1,\dots,f_n)=R^{hs}; $$ par hypothèse les équations $g_1=\cdots=g_n=0$ ont une solution
dans $\widehat{B}$. Comme $B$ est l'hensélisé du localisé d'un schéma de type fini sur un trait complet donc excellent, on peut utiliser 
le théorème d'approximation de M.~Artin. Ainsi, il existe des $x_i$, $1\leq i \leq n+r$, congrus aux $t_i$ modulo $\MM_X^2$ tels que $\mathbf{g}(\mathbf{x})=0$ 
dans $B$. On peut donc définir un $A$-morphisme $\varphi: {A[\mathbf{t}]_{\MM}}^{hs}/(\mathbf{g})\ra B$, par $t_i\mapsto x_i$. Le morphisme $\widehat{\varphi}$ 
induit sur les complétés est un isomorphisme ; le morphisme $\varphi$ est donc étale et, finalement, un isomorphisme.

\subsection{} 
Soient $X,n$, et $\lambda_{X,x}=\lambda$ comme précédemment. Le problème étant local au voisinage de $x$ sur $X$ pour la topologie de Zariski, 
on peut supposer qu'il existe un entier $r$ et
un $S$-module $\mc{E}\iso \OO_S^{n+r}$ tel que $X$ soit isomorphe à un sous-schéma fermé de 
$\mathbf{V}(\mc{E})$, défini par $S$-morphisme $\mathbf{f}:\OO_S^r\ra \mathbf{S}(\mc{E})$, satisfaisant les hypothèses 
de~\ref{not3}. D'après~\ref{deform}, on peut supposer que $\mathbf{f}$ est à valeur dans 
$$\mathbf{S}(\mc{E})_{\leq \lambda}=\OO_S\oplus\cdots\oplus \mathbf{S}^{\lambda}(\mc{E}).$$ Nous noterons $\sur{\mathbf{f}}=\mathbf{f}\otimes_A k$, 
et plus généralement par une barre $-$ toute réduction dans $k$. Remarquons qu'il est important de ne pas choisir immédiatement 
d'isomorphisme $\Gamma(S,\mc{E})\iso A^{n+r}$. Cela simplifie les calculs qui vont suivre. 
(Je dois cette remarque à Luc Illusie.)

Pour tout morphisme $\mathbf{a}:\OO_S^r\ra \mathbf{S}(\mc{E})$, et pour tout entier $i\in \NN$, 
notons $\mathbf{a}^{[i]}$ la composante homogène de degré $i$ de $\mathbf{a}$ et $\mathbf{a}^{[{\leq i}]}=\mathbf{a}^{[0]}+\cdots+\mathbf{a}^{[i]}$ 
sa partie de degré inférieur à $i$. Posons $\widetilde{\mc{E}}=\OO_S t_0 \oplus \mc{E}$, et $\PP_{\mc{E}}=\PP(\widetilde{\mc{E}})$. 
Si $\mathbf{a}$ est un morphisme $\OO_S^r\ra \mathbf{S}(\mc{E})_{\leq \lambda+2}$, notons 
$\widetilde{\mathbf{a}}=\mathbf{a}^{[\lambda+2]}+t_0 \mathbf{a}^{[\lambda]+1}+\cdots+t_0^{\lambda+2}\mathbf{a}^{[0]}$.
Enfin, notons $\varphi_\mathbf{a}:Y(\mathbf{a})\ra S$ le $S$--schéma projectif
$V(\widetilde{\mathbf{a}})\hra \PP_{\mc{E}}$, et $y=(1,0_{\mc{E}})\in (\PP_{\mc{E}})_s$ : c'est
l'image de $x$ par l'immersion composée $X\hra \mathbf{V}(\mc{E})\hra {\PP_{\mc{E}}}$.
Remarquons qu'en vertu de~\ref{deform}, si $\mathbf{a}:\OO_S^r\ra \mathbf{S}(\mc{E})_{\leq \lambda+2}$ a pour $\lambda$--tronqué $\mathbf{f}$,
les hensélisés stricts $Y(\mathbf{a})_{(y)}$ et $X_{(x)}$ sont automatiquement $S$-isomorphes.
On cherche $\mathbf{a}$ tel que $Y(\mathbf{a})$ satisfasse les autres conditions de~\ref{deformproj}, c'est-à-dire les hypothèses de 
lissité et de dimension relative hors de $y$. Il suffit de les vérifier en les point fermés de la fibre spéciale (privée de $y$).
En effet, si elles sont satisfaites en ces points, le schéma $Y(\mathbf{a})$ sera régulier en tous les points fermés de $Y(\mathbf{a})_s$. 
Comme le lieu $\text{reg}(Y(\mathbf{a}))$ des  points réguliers
est ouvert (cf. ÉGA {\sc iv}.6.12.6), cet ouvert contient nécessairement toute la fibre spéciale et, par propreté, on a l'égalité 
$\text{reg}(Y(\mathbf{a}))=Y(\mathbf{a})$.
De même, le morphisme $Y(\mathbf{a})\ra S$ est aussi plat car il est plat en tous les points
fermés de la fibre spéciale et son lieu de platitude est ouvert (cf. ÉGA, {\sc iv}.11.1.1). 
La fibre générique est lisse : tout point $y_{\eta}$ de $Y_{\eta}$ est générisation d'un point $y_s$ de $Y_s$.
Si $y_s$ est différent de $y$, la lissité est évidente par hypothèse tandis que si $y_s=y$ cela résulte du fait 
que l'on a supposé $\SP(\OO_{X,x})-\{x\}$ essentiellement lisse sur~$S$.
Finalement, pour démontrer~\ref{deformproj}, il nous suffit de démontrer la proposition suivante :

\begin{prp}\label{prpfin}
Il existe un morphisme $\mathbf{a}:\OO_S^r\ra \mathbf{S}(\mc{E})_{\leq \lambda+2}$ tel que 
$\mathbf{a}^{[\leq \lambda]}=\mathbf{f}$ et  $Y(\mathbf{a})\times_S s$ soit lisse de dimension
relative $n$ hors de $y$.
\end{prp}

Considérons le $S$--schéma $T=\mathbf{V}\big(\sous{\Hom}(\OO_S^r,S^{\lambda+1}(\mc{E}) \oplus S^{\lambda+2}(\mc{E}))^{\vee}\big)$, paramétrant les morphismes
$\mathbf{a}:\OO_S^r\ra \mathbf{S}(\mc{E})_{\leq \lambda+2}$ tels que $\mathbf{a}^{[\leq \lambda]}=\mathbf{f}$. Enfin, considérons
la variété d'incidence $$M=\{(z,\mathbf{a})\ | \ z\in Y(\mathbf{a}),\text{\ non lisse de dim. rel.}\ n \text{\ en}\ z\}\hra ({\PP_{\mc{E}}}-\{y\})\times_S T,$$
sa structure de schéma est précisée plus loin, à l'aide du critère Jacobien.
Elle est naturellement munie de deux projections $p_1:M\ra {\PP_{\mc{E}}}$ et $p_2:M\ra T$.
La proposition précédente est une conséquence immédiate du lemme suivant :
\begin{lmm}
Pour tout $z\in (\PP_{\mc{E}})_s$, $z\neq y$, on a $\dim\,p_1^{-1}(z)=\dim\,T_s-n-r-1$.
\end{lmm}
En effet, on aura alors $\dim\,M_s\leq \dim\,T_s-1$ et finalement $\dim\,p_2(M_s)<\dim\,T_s$. En particulier, on
aura $p_2(M(k))\neq T(k)$, ce qui démontre~\ref{prpfin} en relevant arbitrairement un élément de $T(k)\diagdown p_2(M(k))$.

\begin{proof}
Soit $z$ comme dans l'énoncé. Pour démontrer le lemme, on peut choisir un isomorphisme $\widetilde{\mc{E}}\iso \OO_S t_0\oplus\cdots\oplus \OO_S t_{n+r}$
tel que $z\in (\PP_{\mc{E}})_s\iso \PP^{n+r}_s$ (resp. $y$) soit de coordonnées $(0,0,\dots,0,1)$ (resp. $(1,0,0,\dots,0)$). 
Cela résulte du fait que $z$ est supposé différent de $y$.
Il est équivalent de se donner $\mathbf{a}$ dans $T(S)$ (resp. dans $T(s)$) et, pour chaque $i\in \{1,\dots,r\}$ et chaque suite (finie) 
$\alpha\in \NN^{n+r}$ de somme $|\alpha|\in \{\lambda+1,\lambda+2\}$, un élément $c_{\mathbf{a},i,\alpha}$ dans $A$ (resp. dans $k$). 
À une telle famille de coefficients, on associe $\widetilde{\mathbf{a}}=(\widetilde{\mathbf{a}}_1,\dots,
\widetilde{\mathbf{a}}_r)$,
où 
$$
\widetilde{\mathbf{a}}_i=\underbrace{t_0^{\lambda+2}\mathbf{f}_i(\frac{t_1}{t_{0}},\dots,\frac{t_{n+r}}{t_0})}_{\mathbf{g}_i(t_0,\dots,t_{n+r})}+
\sum_{|\alpha|\in \{\lambda+1,\lambda+2\}}
c_{\mathbf{a},i,\alpha} \underbrace{t_0^{\lambda+2-|\alpha|}t_1^{\alpha_1}\cdots t_{n+r}^{\alpha_{n+r}}}_{m_{\alpha}(\mathbf{t})}$$
(resp. la même expression avec $\sur{\mathbf{f}}_i$ à la place de $\mathbf{f}_i$).
Considérons maintenant $\mathbf{a}\in p_1^{-1}(z)$. Nécessairement, pour chaque $i\in \{1,\dots,r\}$, le seul monôme
$m_{\alpha}(\mathbf{t})$ qui ne soit pas nul évalué en $z$ est celui pour lequel $\alpha=(0,0,\dots,0,\lambda+2)\fed \beta$.
La condition $z\in Y(\mathbf{a})$ s'écrit donc $\sur{\mathbf{g}}_i(z_0,z_1,\dots,z_{n+r-1},1)+c_{\mathbf{a},i,\beta}=0$, 
pour chaque $i\in \{1,\dots,r\}$.
Plaçons nous dans l'ouvert affine $t_{n+r}\neq 0$ de $\PP^{n+r}_s$ ; 
il contient $z$ par hypothèse. Notons $t_0',t_1',\cdots,t_{n+r-1}'$ les coordonnées affines déduites
de $t_0,\dots,t_{n+r}$ (ainsi, $t_i'=\frac{t_i}{t_{n+r}}$).
Pour $i\in \{1,\dots,r\}$, calculons la dérivée partielle $\displaystyle \frac{\partial \widetilde{\mathbf{a}}_i}{\partial t'_0}$ en $z$ dans ces coordonnées.
Elle vaut :
$$
\underbrace{\frac{\partial \sur{\mathbf{g}}_i}{\partial t'_0}\overbrace{(0,\dots,0)}^{n+r\ \text{zéros}}}_{\chi_{i,0}}+c_{\mathbf{a},i,\gamma},
$$
où $\gamma\in \NN^{n+r}$ est la suite $(0,\dots,0,\lambda+1)$. En effet, pour que la dérivée partielle par rapport à $t'_0$ du monôme $t_0^{\prime\lambda+2-|\alpha|} 
t_1^{\prime\alpha_1}\cdots t_{n+r-1}^{\prime\alpha_{n+r-1}}$ soit non nulle évaluée en $(0,\dots,0)$, il faut que $\alpha_1=\cdots=\alpha_{n+r-1}=0$
et $|\alpha|=\lambda+1$. Pour $\alpha=\gamma$ cette dérivée partielle vaut $1$. 
Si maintenant $j$ est un indice dans $\{1,\dots,n+r-1\}$, on a de même :
$$
\frac{\partial \widetilde{\mathbf{a}}_i}{\partial t'_j}(0,\dots,0)
=\underbrace{\frac{\partial \sur{\mathbf{g}}_i}{\partial t'_j}(0,\dots,0)}_{\chi_{i,j}}
+c_{\mathbf{a},i,\gamma(j)}
,$$
où $\gamma(j)\in \NN^{n+r}$ est définie par $\gamma(j)_j=1$, $\gamma(j)_u=0$ pour $u\notin \{j,n+r\}$, et $|\gamma(j)|=\lambda+2$.
On peut noter qu'avec ces conventions, on a l'égalité $\chi_{i,j}=0$ pour tous les couples $(i,j)$ considérés. Cela résulte du fait
que $t_0^{\prime 2}$ divise tous les $\sur{\mathbf{g}_i}$.
Rappelons que si $\mathbf{a}\in T(s)$, et $z\in Y(\mathbf{a})(s)$, 
le morphisme $\varphi_{\mathbf{a}}$ est lisse de dimension relative $n$ en $z$ si, et seulement si,
le rang en $z$ d'une matrice Jacobienne des $\widetilde{\mathbf{a}}_i$ (dans une carte affine quelconque contenant $z$) est égal à $r$.
Finalement, la condition $\mathbf{a}\in M_z$ est définie par l'intersection du sous-espace affine $L_z$ de codimension $r$ de $T_s$
d'équations $\sur{\mathbf{g}}_i(0,0,\dots,0,1)+c_{\mathbf{a},i,\beta}=0$ ($i\in \{1,\dots,r\})$
et du sous-schéma de $T_s$ d'équations les mineurs $r\times r$ de la matrice 
$$
\begin{pmatrix}
 c_{\mathbf{a},1,\gamma} & \cdots & c_{\mathbf{a},r,\gamma}\\
c_{\mathbf{a},1,\gamma(1)} & \cdots & c_{\mathbf{a},r,\gamma(1)}\\
   \vdots & \ddots & \vdots \\
c_{\mathbf{a},1,\gamma({n+r-1})} & \cdots & c_{\mathbf{a},r,\gamma({n+r-1})}\\
\end{pmatrix}.
$$
Comme les suites $\gamma,\gamma(i)$, $i\in \{1,\dots,n+r-1\}$ sont distinctes, le sous-schéma $\mathrm{Min}$ défini par l'annulation de ces mineurs est
de codimension $n+1$ dans l'espace affine $T_s$ (cf. par exemple~\cite{artintata}). De plus comme ces suites sont différentes de $\beta$,
l'intersection de $\mathrm{Min}$ avec le sous-espace affine $L_z$ est transverse. Ainsi, $\mathrm{codim}_{T_s}(M_z)=r+(n+1)$, d'où le résultat.
\end{proof}

\subsection{Le cas de la dimension relative $1$}\label{Michel}Voici l'argument de Michel Raynaud qui permet de démontrer directement~\ref{deformproj}
dans le cas des courbes. Supposons $X/S$ affine (et $S$ complet). Notons $Y=X_s$ la fibre spéciale. Comme c'est une courbe, il existe une compactification $Z$
de $Y$, projective, et lisse hors de $x$. Le $s$-schéma $Z-\{x\}$ est affine et lisse donc (cf.~\sga{1}{iii}{6.8}) il existe un $S$-schéma formel affine et lisse
$\mc{T}$ dont $Z-\{x\}$ est la fibre spéciale. L'affine $Y-\{x\}$ admet un unique relèvement formel $\mc{U}$ sur $S$. Il est  
naturellement muni d'immersions ouvertes $\mc{U}\hra \mc{T}$ et $\mc{U}\hra \widehat{X}$. On peut donc recoller $\mc{T}$ et $\widehat{X}$ le long de $\mc{U}$ :
le schéma formel $\mc{T}\coprod_{\mc{U}} \widehat{X}$ est une déformation plate de $Z$, donc propre sur $S$ 
et s'algébrise (cf.~\sga{1}{iii}{7.2}) en un schéma $X'$ sur $S$
propre et lisse hors de $x$ sur $S$, qui est formellement, donc localement pour la topologie étale (cf.~\ref{algebris}), isomorphe à $X$ en $x$.
\section{Démonstration du théorème~\ref{theoreme}}
Commençons par remarquer que pour démontrer la conjecture~\ref{conjDM}, on peut supposer $S$ complet.
Il est bien connu que le terme étale est invariant par une telle extension $\hat{S}\ra S$ (cf.~\sga{$\mathbf{4\frac{1}{2}}$}{Th. Finitude}{3.7}). 
L'égalité $\mu(X/S,\, x)=\mu(X_{\hat{S}}/\hat{S},\,x)$
résulte de l'isomorphisme $\OO_{Z_{\hat{S}}}=\OO_{Z_S}\otimes_{\OO_S} \OO_{\hat{S}}$, dans les notations de~\ref{t1local}.
Ceci étant, on peut supposer d'après~\ref{deformproj}, que $X/S$ est propre car les deux termes de l'égalité à démontrer ne dépendent que de
l'hensélisé (strict) en $x$. 
D'un côté on a inconditionnellement, 
$$\chi(X,\mathsf{L}\Lambda^{n+1}\Om^1_{X/S})\sr{\ref{globalisation}}{=}\chi(X,\Om^{n+1}_{X/S})\sr{\text{\ref{t1local}}}{=}\mu(X/S,\,x),$$
tandis que la conjecture de Bloch prédit que 
$$\chi(X,\mathsf{L}\Lambda^{n+1}\Om^1_{X/S})=(-1)^n\mathrm{Art}(X/S)=(-1)^n\text{dim\ tot}\, \ce(\FF_{\ell}).\text{\ \sc c.q.f.d.}$$ 
\bibliography{bib}
\bibliographystyle{smfalpha}

\end{document}

%% file: stylemoi2
\def\doute#1{}

\usepackage{amsmath,a4wide}
\usepackage{amssymb}
\usepackage{mathrsfs}
\usepackage{inputenc,xspace}
\usepackage{euscript}
\usepackage[all]{xy}
\usepackage{mathrsfs}
\usepackage{french}
\usepackage{smfthm}

\SwapTheoremNumbers

\makeatletter
\renewcommand\theequation{\thesubsection.\arabic{equation}}
\@addtoreset{equation}{subsection}
\makeatother

\newtheorem{prp}[subsection]{Proposition} 
\newtheorem{thm}[subsection]{Théorème}
\newtheorem{lmm}[subsection]{Lemme}
\newtheorem{cnj}[subsection]{Conjecture}
\newtheorem{crl}[subsection]{Corollaire}

\theoremstyle{definition}

\def\sga#1#2#3{[{\bf SGA~#1}~{\sc #2}~#3]}
\newcommand{\beqt}{\begin{equation}}
\newcommand{\eeqt}{\end{equation}}

\newcommand{\Hom}{\mathsf{Hom}}

\newcommand{\K}{{\mathsf{K}}}

\newcommand{\OO}{\mathscr{O}} 

\newcommand{\SP}{\mathrm{Spec}}
\newcommand{\aff}{\mathbf{A}}
\newcommand{\etb}{\bar{\eta}}

\newcommand{\red}{\text{réd}}
\newcommand{\Om}{{\Omega}}

\newcommand{\ce}{\Phi}
\newcommand{\RG}{\mathrm{R\Gamma}}

\newcommand{\HH}{\mathrm{H}}

\newcommand{\swan}{\mathrm{Swan}}

\newcommand{\gp}{\pi_1}     

\newcommand{\Ext}{{\mathsf{Ext}}}

\newcommand{\kos}{{\mathsf{Kos}^{\llcorner}}}
\newcommand{\kosdual}{{\mathsf{Kos}^{\wedge}}}

\newcommand{\ZZ}{\mathbf{Z}}

\newcommand{\NN}{\mathbf{N}}

\newcommand{\PP}{\mathbf{P}}
\newcommand{\FF}{\mathbf{F}}

\newcommand{\MM}{\mathfrak{m}}

\newcommand{\ra}{\rightarrow}

\newcommand{\hra}{\hookrightarrow }
\newcommand{\sr}{\stackrel}

\newcommand{\rrraxy}{\ar@<1ex>[r] \ar@<-1ex>[r] \ar[r] }
\newcommand{\iso}{\stackrel{\sim}{\ra}}

\newcommand{\fed}{\stackrel{\text{\rm déf}}{=}}

\def\]{\textup{\mbox{]\hspace{-.15em}]}}}
\def\[{\textup{\mbox{[\hspace{-.15em}[}}}
\def\mc{\mathscr}

\def\-l{\ \vspace{-4mm}}

\def\L#1{\mathsf{L}_{#1}}

\def\sur{\overline}
\def\sous{\underline}

\author{Fabrice Orgogozo}
\address{École Normale Supérieure, DMA \\ 45, Rue d'Ulm \\F-75230 Paris,  Cedex 05 \\ France}
\email{fabrice.orgogozo@ens.fr}

\def\nL2#1{{\vert \vert #1 \vert \vert }_{\mathrm{L}^2}}
\def\n2t#1{w_{\mathrm{L}^2,\iota}(#1)}